\documentclass{article}
\advance \topmargin by -\headheight
\advance \topmargin by -\headsep
     
\evensidemargin \oddsidemargin
\newcommand\be{\begin{equation}}
\newcommand\ee{\end{equation}}

\begin{document}

\title{A recursion equation for prime numbers}         
\author{Joseph B.\ Keller\\Department of Mathematics\\Stanford University\\ Stanford, CA  94305\\Email: keller@math.stanford.edu}        

\date{June 24, 2007}
\maketitle

\begin{abstract}
It is shown that the first $n$ prime numbers $p_1,\ldots,p_n$ determine the next one by the recursion equation
$$
p_{n+1} =\lim\limits_{s\to +\infty} \left[\prod\limits^n_{k=1} \left( 1-\frac{1}{p^s_k}\right)\;\sum\limits^\infty_{j=1} \;\frac{1}{j^s} -1\right]^{-1/s}.
$$
The upper limit on the sum can be replaced by  $2p_n -1$, and the result still holds.
\end{abstract}

\vspace{.5in}

The first $n$ prime numbers $p_1, p_2,\cdots,p_n$ determine the next one by the  recursion equation:
\be
p_{n+1} =\lim\limits_{s\to +\infty} \left[\prod\limits^n_{k=1} \left( 1-\frac{1}{p^s_k}\right)\;\sum\limits^\infty_{j=1} \;\frac{1}{j^s} -1\right]^{-1/s}.
\label{1}
\ee

To prove (\ref{1}) we replace the sum in it, which is the zeta function $\zeta(s)$, by the Euler product\cite{Edwards}, to get
\be
\prod\limits^n_{k=1} \left( 1-\frac{1}{p^s_k}\right) \zeta(s) =\prod\limits^n_{k=1} \left( 1-\frac{1}{p^s_k} \right) \prod\limits^\infty_{j=1} \left( 1-\frac{1}{p^s_j} \right)^{-1} =\prod\limits^\infty_{j=n+1} \left( 1-\frac{1}{p^s_j} \right)^{-1}.
\label{2}
\ee
As $s\to +\infty$, the final product in (\ref{2}) is asymptotic to $1+p^{-s}_{n+1}$.  When $1$ is subtracted from it, the difference is asymptotic to $p^{-s}_{n+1}$.  Taking the $-1/s$ power of this yields $p_{n+1}$, which proves (\ref{1}).

We can obtain another recursion equation by first factoring the product out of the bracketed expression in (\ref{1}).  Then as $s\to \infty$, the product tends to $1$, and so does the $-1/s$ power of it.  Thus we obtain the recursion equation
\be
p_{n+1} =\lim\limits_{s\to +\infty} \; \left[ \sum\limits^\infty_{j=1} \; \frac{1}{j^s} - \prod\limits^n_{k=1} \left( 1-\frac{1}{p^s_k} \right)^{-1} \right]^{-1/s}.
\label{3}
\ee
Since $p_{n+1} < 2p_n$, (\ref{3}) remains valid when the sum is terminated at $j=2p_n -1$, so it becomes a finite sum:
\be
p_{n+1} =\lim\limits_{s\to \infty} \left[ \sum\limits^{2p_n-1}_{j=1} \: \frac{1}{j^s} -\prod\limits^n_{k=1} \; \left( 1-\frac{1}{p^s_k} \right)^{-1} \right]^{-1/s} .
\label{4}
\ee

To show that (\ref{1}) also holds with the sum terminated at $2p_n-1$, we just factor the product out of the bracketed expression in (\ref{4}).  The $-1/s$ power of  the product tends to $1$ as $s$ tends to $+\infty$, so the result follows.


\begin{thebibliography}{99}
\bibitem{Edwards} Edwards, Harold.  {\em Riemann's Zeta Function}, Academic Press, New York, 1974.
\end{thebibliography}
\end{document}